\documentclass{article}
\usepackage{fullpage,amssymb}
\newenvironment{@abssec}[1]{%
     \if@twocolumn
       \section*{#1}%
     \else
       \vspace{.05in}\footnotesize
       \parindent .2in
         {\bfseries #1. }\ignorespaces
     \fi}
     {\if@twocolumn\else\par\vspace{.1in}\fi}
\newenvironment{keywords}{\begin{@abssec}{Key words}}{\end{@abssec}}
\newenvironment{AMS}{\begin{@abssec}{AMS subject classification}}{\end{@abssec}}
\newenvironment{proof}{\noindent{\bf Proof.\/}}{\hfill $\Box$}

\newtheorem{theorem}{Theorem}[section]

\newtheorem{corollary}[theorem]{Corollary}

\begin{document}

\newfont{\bb}{msbm10}
\def\Bbb#1{\mbox{\bb #1}}
\def\C{{\Bbb C}}
\def\N{{\Bbb N}}
\def\order{\mathop{\rm order}\nolimits}
\def\belowrightarrow#1{{{{}\over\ #1\ }\kern-1.1em\to}}
\def\eqbd{\mathop{{:}{=}}\nolimits}
\def\bdeq{\mathop{{=}{:}}\nolimits}

\title{On convergence of infinite matrix products}
\author{Olga Holtz  \\
Department of Computer Sciences, \\
 University of Wisconsin, Madison, Wisconsin 53706 U.S.A. \\ 
holtz@cs.wisc.edu\thanks{This work was supported in
part by the Clay Mathematics Institute.}}
\date{}
\maketitle

\begin{abstract} A necessary and sufficient condition for the convergence of
an infinite right product of matrices of the form
$$ A \eqbd \left[ \begin{array}{cc} I & B \\ 0 & C \end{array} \right], $$ 
with (uniformly) contracting submatrices $C$, is proven.
\end{abstract}

\begin{keywords} Infinite matrix products, RCP sets.
\end{keywords}

\begin{AMS} 15A60, 15A99
\end{AMS}

\section{Introduction}

Consider the set of all matrices in $\C^{d\times d}$ 
of the form   \begin{equation} 
 A\eqbd \left[ \begin{array}{cc} I_s & B \\ 0 & C 
\end{array} \right],  \label{triag} 
\end{equation}
where $I_s$ denotes the identity matrix of order $s<d$.

Matrices~(\ref{triag}) are known (e.g.,~\cite{DL}) to form an LCP set whenever
the submatrices $B$ are uniformly bounded and the submatrices $C$ 
are uniformly contracting, that is,  satisfy the condition $\|C\|\leq r$ 
for some fixed matrix (i.e., submultiplicative) norm  $\| \cdot \|$ on 
$\C^{(d-s)\times (d-s)}$ and  some constant $r<1$.  To recall, a set 
$\Sigma$ has the LCP (RCP) property  if all left (right) infinite 
products formed  from matrices in $\Sigma$  are convergent. 

Matrices of the form~(\ref{triag}), with uniformly bounded submatrices
$B$ and uniformly contracting submatrices $C$, do not necessarily form an
RCP set. (They do form such a set if and only if they satisfy a very 
stringent condition given in Corollary~\ref{cor} below.) However, there exists
a simple criterion that can be used to check whether a {\em particular\/} 
right  infinite product formed from such matrices converges.

\section{A convergence test}

\begin{theorem} \label{thm} Let $(A_n)_{n\in \N}$ be a sequence of matrices
of the form~(\ref{triag})  and let $$ \|C_n\|\leq r<1 \qquad \mbox{\rm for all}
\quad n\in \N$$ for some matrix norm $\| \cdot \|$. The  sequence
$( P_n\eqbd A_1A_2\cdots A_n )$ converges if and only if
so does the sequence $(B_n(I-C_n)^{-1})$. In this event,
$$ \lim_{n \to \infty} P_n=\left[ \begin{array}{cc}
                                I & \lim_{n\to \infty} B_n(I-C_n)^{-1} \\
                                0 & 0 \end{array}\right]. $$ 
\end{theorem}

\begin{proof}
To prove  the necessity, partition $P_n$  conformably  with $A_n$.  Then
$$ P_n=\left[  \begin{array}{cc} I & X_n  \\ 0 & C_1 C_2 \cdots C_n
 \end{array} \right] \qquad \mbox{where} \qquad X_n\eqbd  
\sum_{i=0}^n B_{n-i} (C_{n+1-i} C_{n+2-i} \cdots C_n).$$ 
If $(P_n)$ converges, then $\lim_{n\to \infty}(X_n-X_{n-1})=0$.
Also, $\|(I-C_n)^{-1}\|\leq 1/(1-r)$ for all $n\in \N$. 
But $X_n=B_n+X_{n-1}C_n$, so 
$$ B_n(I-C_n)^{-1}-X_{n-1}=(X_n-X_{n-1})(I-C_n)^{-1} \belowrightarrow
{n\to \infty} 0.$$
Hence $$\lim_{n\to \infty} B_n(I-C_n)^{-1}=\lim_{n\to \infty} X_n.$$

Now prove the sufficiency. Without loss of generality one can assume
that $s=d-s$. Indeed, simply replace each $A_n$ by 
$$ \widetilde{A_n}\eqbd \left[ \begin{array}{cc} I_{\max\{s,d-s\}} & 
\widetilde{B_n}   \\ 0 & \widetilde{C_n} \end{array} \right] $$
where
\begin{eqnarray*}
 \widetilde{B_n} & \eqbd & \cases{ \left[\begin{array}{cc} 
B_n & 0_{s\times (2s-d)} \end{array} \right] & if $s\geq d-s$ \cr
\left[ \begin{array}{c} B_n \\ 0_{(d-2s)\times (d-s)} \end{array}
\right] & if $s< d-s$  }, \\
 \widetilde{C_n} &\eqbd& \cases{ \left[\begin{array}{cc} 
C_n & 0_{(d-s)\times (2s-d)} \\
0_{(2s-d)\times(d-s)} & 0_{2s-d} \end{array} \right] & if $s\geq d-s$ \cr
C_n & if $s< d-s$  }. \end{eqnarray*} 
Then the matrices $\widetilde{A_n}$ satisfy all the assumptions of the theorem
and the sequence $(B_n(I-C_n)^{-1})$ (the product $P_n$) converges iff so 
does the sequence $(\widetilde{B_n}(I-\widetilde{C_n})^{-1})$  
(the product $\widetilde{P_n}$). 

Thus, assume that $s=d-s$. Note that if the sequence 
$(B_n(I-C_n)^{-1})$ converges, then the sequence $(B_n)$ is bounded,
since $\|I-C_n\|\leq 1+r$ for all $n$.
Now, let  \begin{eqnarray*}
D_n & \eqbd & X_n - B_n  (I-C_n)^{-1} \\  
Y_n & \eqbd & B_{n+1} (I-C_{n+1})^{-1} -B_n(I-C_n)^{-1}  
\end{eqnarray*}
for all $ n\in \N$. Then  
\begin{equation} D_{n+1}=(D_n-Y_n)C_{n+1}, \label{iden} \end{equation}
hence
$$ \| D_{n+1} \| \leq (\|D_n\|+\|Y_n\|) \|C_{n+1}\| 
\leq  (\|D_n\|+\|Y_n\|) r. $$
Repeated use of this inequality gives
$$  \|D_n\|\leq \sum_{i=1}^{n-1}  \|Y_{n-i}\| r^i. $$
This implies, in particular, that $$S\eqbd \lim\sup_{n\to \infty} \|D_n\|<
\infty.$$ Since $\lim_{n\to \infty} Y_n=0$, the identity~(\ref{iden}) and 
the upper  bound on $\|C_n\|$ imply that $S\leq rS$, therefore $S=0$, that is, 
$$\lim_{n\to \infty} D_n=0.$$ \end{proof}

The obtained criterion of convergence can be used to make two more
observations in the same spirit.

\begin{corollary} Let $(A_n)_{n\in \N}$ be a sequence of matrices
of the form~(\ref{triag}) such that the sequence $(C_n)$ converges to
a matrix $C$ with spectral radius smaller than $1$. Then the  sequence
$( P_n\eqbd A_1A_2\cdots A_n )$ converges if and only if
so does the sequence $(B_n)$. In this event,
$$ \lim_{n \to \infty} P_n=\left[ \begin{array}{cc}
                                I & \lim_{n\to \infty} B_n  (I-C)^{-1} \\
                                0 & 0 \end{array}\right]. $$ 
\end{corollary}

\begin{proof}
If $\varrho(C)<1$, then there exists a matrix norm $\| \cdot \|$ 
on $\C^{(d-s)\times (d-s)}$ such that $\| C\|<1$ 
(e.g.,~[p.297, Lemma 5.6.10]\cite{HJ}).
So, $\|C_n\|\leq r$ for all $n\geq N$ for some $r<1$ and
some $N\in \N$, so the assumption of the theorem is then satisfied.
The product $P_n$ converges whenever so does the product $A_{N} A_{N+1}
 \cdots$, so $(P_n)$ has a limit whenever $(B_n)$ has one. 
By the same reason, the sequence $((I-C_n)^{-1})_{n=N}^\infty$
is bounded, so the necessity argument from the proof of the Theorem
shows that the convergence of $(B_n)$ is also necessary. \end{proof}

\begin{corollary} \label{cor} A set $\Sigma$ consisting of matrices
of the form~(\ref{triag},) with uniformly contracting submatrices 
$C$, is an RCP set  if and only if 
\begin{equation} 
 B_1(I-C_1)^{-1}=B_2 (I - C_2)^{-1} \qquad 
\mbox{\rm for all}\quad  A_1, \; A_2 \in \Sigma, \label{crit} \end{equation}
where $$ A_i = \left[ \begin{array}{cc} I & B_i \\ 0 & C_i 
\end{array} \right], \qquad i=1,2.$$ 
\end{corollary}

\begin{proof} Given $A_1$, $A_2\in \Sigma$, apply Theorem~\ref{thm} to the 
product  $A_1 A_2 A_1 A_2 \cdots$  to see that the 
condition~(\ref{crit}) is necessary and sufficient for the convergence
of such a product. But if it is satisfied for all pairs of matrices
from $\Sigma$, then it is sufficient for the convergence of any
right product of matrices from $\Sigma$. \end{proof}

\section*{Acknowledgements} I am grateful to Professor Hans Schneider 
for his critical reading of the manuscript and to the referee for 
valuable suggestions.

\end{document}